\documentclass[10pt,a4paper,reqno]{amsart}
\usepackage{amssymb,amsmath,amsthm}
\usepackage[T1]{fontenc}
\usepackage{graphicx}
\usepackage{caption}

\theoremstyle{plain}

\title{Other trigonometric proofs of Pythagoras theorem}
\author{Nuno Luzia}
\address{Universidade Federal do Rio de Janeiro, Instituto de Matem\'atica \\ Rio de Janeiro 21941-909, Brazil}  
\email{nuno@im.ufrj.br} 

\begin{document}
\maketitle 
\begin{abstract}
Only very recently a trigonometric proof of the Pythagoras theorem was given by Zimba \cite{1}, many authors thought this was not possible. In this note we give other trigonometric proofs of Pythagoras theorem by establishing, geometrically, the half-angle formula $\cos\theta=1-2\sin^2 \frac{\theta}{2}$.
\end{abstract}

\section{Pythagoras theorem}
Pythagoras theorem is ilustrated in figure below.

\begin{figure}[h]
\begin{center}
  \includegraphics[scale=0.6]{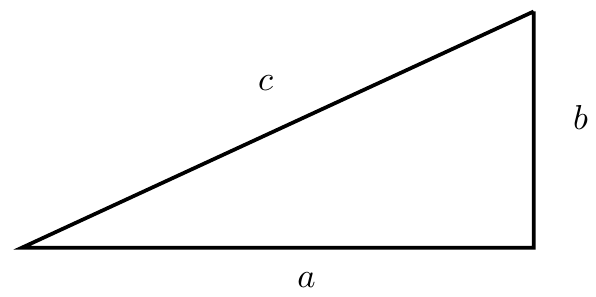}
 \caption*{$a^2+b^2=c^2$}
\end{center}
\end{figure}
To prove this theorem it is enough to consider the particular case $c=1$, for then the general case follows by rescaling. Define the trigonometric functions $\sin \theta$ and
$\cos \theta$, for $0<\theta<90^\circ$, as usual (see Figure 1).

\begin{figure}[h]
\begin{center}
  \includegraphics[scale=0.6]{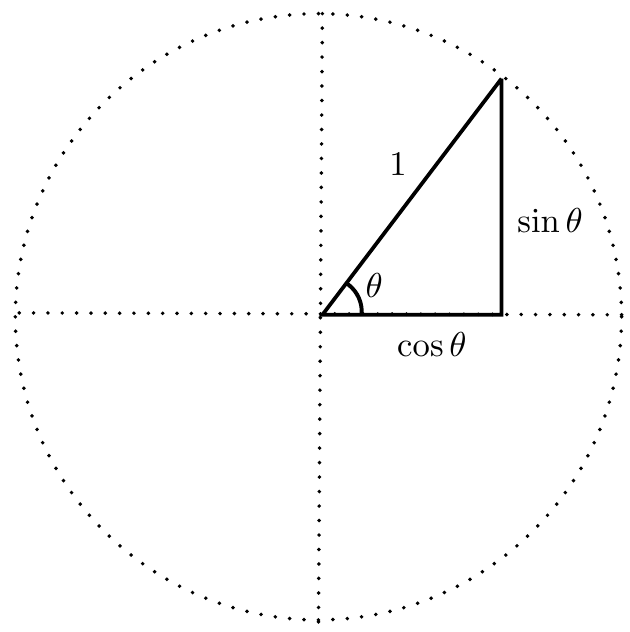}
 \caption{Definition of $\sin \theta$ and $\cos \theta$}
\end{center}
\end{figure}
Then Pythagoras theorem is equivalent to
\begin{equation}\label{fund}
   \cos^2 \theta + \sin^2 \theta=1,
\end{equation}
for every $0<\theta<90^\circ$.

The subtracting formulas
\begin{align*}
&\cos(\alpha-\beta)=\cos\alpha\cos\beta+\sin\alpha\sin\beta,\\
&\sin(\alpha-\beta)=\sin\alpha\cos\beta-\cos\alpha\sin\beta,
\end{align*}
valid for $0<\alpha,\beta,\alpha-\beta<90^\circ$, can be proved geometrically without using Pythagoras theorem (see \cite{1}). Zimba \cite{1} used these formulas to prove Pythagoras theorem.

The addition formula
\[
     \cos(\alpha+\beta)=\cos\alpha\cos\beta-\sin\alpha\sin\beta,
\]
$0<\alpha,\beta,\alpha+\beta<90^\circ$, can, in the same way, also be proved geometrically without using Pythagoras theorem, from which it follows the half-angle formula
\begin{equation}\label{half1}
\cos\theta=\cos^2 \frac{\theta}{2}-\sin^2\frac{\theta}{2}.
\end{equation}
If we use Pythagoras formula (\ref{fund}) in (\ref{half1}) then we get
\begin{equation}\label{half2}
     \cos\theta=1-2\sin^2 \frac{\theta}{2}.
\end{equation}
What if we could prove (\ref{half2}) without using Pythagoras theorem? Then using  (\ref{half1}) and (\ref{half2}) we would get the Pythagoras formula (\ref{fund}) !

\subsection{Proof of half-angle formulas}
First we observe the simple fact that in an isosceles triangle with two equal sides with length 1, forming an angle $\theta$, the length of the other side is $2\sin\frac{\theta}{2}$. This can be seen by dividing the isosceles triangle in two right triangles with an angle $\frac{\theta}{2}$. See Figure 2.

\begin{figure}[h]
\begin{center}
  \includegraphics[scale=0.6]{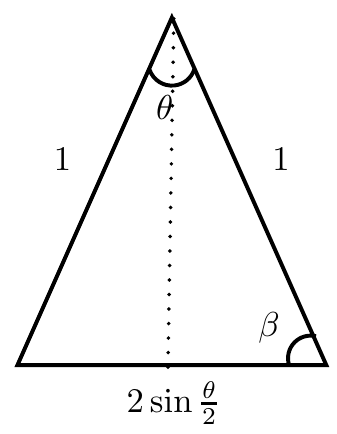}
 \caption{}
\end{center}
\end{figure}
Also the other angles $\beta$ of the isosceles triangle satisfy $\theta+\beta+\beta=180^\circ$ or $\beta=90^\circ-\frac{\theta}{2}$. Of course, by definition, $\cos (90^\circ-\theta)=\sin\theta$ and $\sin (90^\circ-\theta)=\cos\theta$.

Now consider the two right triangles as in Figure 3, which together form an isosceles triangle as in Figure 2.

\begin{figure}[h]
\begin{center}
  \includegraphics[scale=0.7]{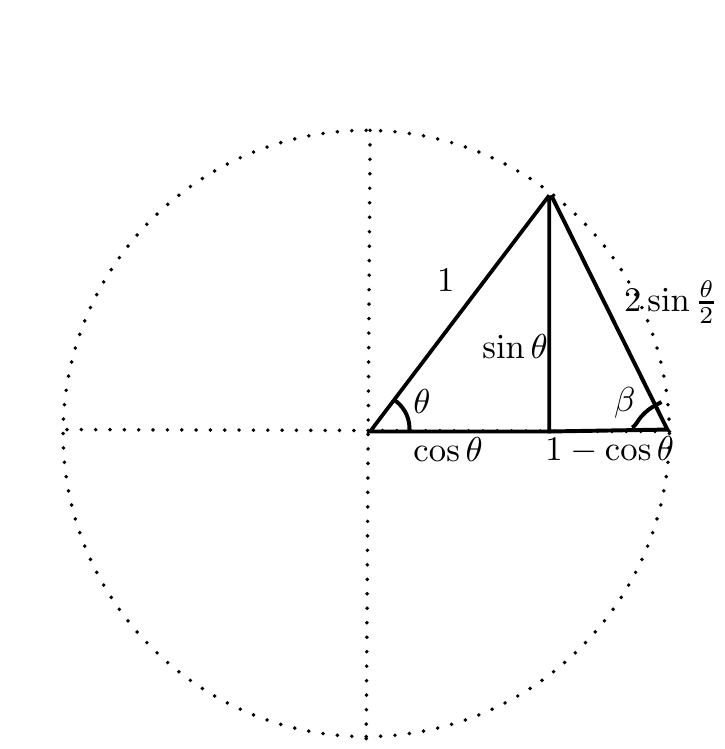}
 \caption{}
\end{center}
\end{figure}

Using the trigonometric relations for the right hand side triangle, we get
\begin{align*}
   1-\cos\theta&=2\sin\frac{\theta}{2}\cos(90^\circ-\frac{\theta}{2}),\\
   \sin\theta&=2\sin\frac{\theta}{2} \sin(90^\circ-\frac{\theta}{2}),
\end{align*}
 or
\begin{align}
   \cos\theta&=1-2\sin^2\frac{\theta}{2}, \label{half3}\\
   \sin\theta&=2\sin\frac{\theta}{2} \cos\frac{\theta}{2}. \label{half4}
\end{align}
 
As mentioned before, (\ref{half3}) together with  (\ref{half1}) imply the Pythagoras formula (\ref{fund}).

\subsection{Another proof}
Here we prove the Pythagoras theorem using formulas (\ref{half3}) and (\ref{half4}) (instead of using formulas (\ref{half3}) and (\ref{half1})). We get
\[
     \cos^2 \theta + \sin^2\theta=1-4\sin^2\frac{\theta}{2} \left( \cos^2\frac{\theta}{2}+\sin^2\frac{\theta}{2}-1\right).
\]
By mathematical induction we get
\begin{equation}\label{ind}
     \cos^2 \theta + \sin^2\theta=1-4^n\left(\prod_{i=1}^n\sin^2\frac{\theta}{2^i}\right) \left( \cos^2\frac{\theta}{2^n}+\sin^2\frac{\theta}{2^n}-1\right),
\end{equation}
for every natural number $n$. Clearly $\left| \cos^2\frac{\theta}{2^n}+\sin^2\frac{\theta}{2^n}-1\right|\le 1$ (actually this quantity can be made arbitrarly small by making $n$ large). Let $\theta$ be measured in radians. It follows from Figure 3 that $\sin\theta<2\sin\frac{\theta}{2}<\theta$. Then
\[
     4^n\prod_{i=1}^n\sin^2\frac{\theta}{2^i}<(4\theta^2)^n 4^{-\sum_{i=1}^n i} < (4\theta^2)^n 2^{-n^2},
\]
which is clearly arbitrarly small by making $n$ large enough. Then (\ref{ind}) implies the Pythagoras formula (\ref{fund}).

\end{document}